\documentclass[sn-mathphys-num]{sn-jnl}

\usepackage{graphicx}%
\usepackage{multirow}%
\usepackage{amsmath,amssymb,amsfonts}%
\usepackage{amsthm}%
\usepackage{mathrsfs}%
\usepackage[title]{appendix}%
\usepackage{xcolor}%
\usepackage{textcomp}%
\usepackage{manyfoot}%
\usepackage{booktabs}%
\usepackage{algorithm}%
\usepackage{algorithmicx}%
\usepackage{algpseudocode}%
\usepackage{listings}%

\theoremstyle{thmstyleone}%
\newtheorem{theorem}{Theorem}
\newtheorem{prop}{Proposition}[section]%
\newtheorem{lem}{Lemma}[section]

\theoremstyle{thmstyletwo}%
\newtheorem{question}{Question}

\theoremstyle{thmstylethree}%

\newcommand{\Aut}{\operatorname{Aut}}

\begin{document}

\title[Article Title]{On combinatorial properties of Gruenberg--Kegel graphs of finite groups}


\author[1]{\fnm{Mingzhu} \sur{Chen}}\email{994194@hainanu.edu.cn}
\equalcont{These authors contributed equally to this work.}

\author[2,3]{\fnm{Ilya} \sur{Gorshkov}}\email{ilygor8@gmail.com}
\equalcont{These authors contributed equally to this work.}

\author[4,5]{\fnm{Natalia~V.} \sur{Maslova}}\email{butterson@mail.ru}
\equalcont{These authors contributed equally to this work.}

\author[6]{\fnm{Nanying} \sur{Yang}}\email{yangny@jiangnan.edu.cn}
\equalcont{These authors contributed equally to this work.}

\affil[1]{\orgdiv{School of Mathematics and Statistics}, \orgname{Hainan University}, \orgaddress{\city{Haikou}, \postcode{570225}, \state{Hainan}, \country{P. R. China}}}

\affil[2]{\orgname{Sobolev Institute of Mathematics SB RAS}, \orgaddress{\city{Novosibirsk}, \postcode{630090}, \country{Russia}}}

\affil[3]{\orgname{Siberian Federal University}, \orgaddress{\city{Krasnoyarsk}, \postcode{660041}, \country{Russia}}}

\affil[4]{\orgname{Krasovskii Institute of Mathematics and Mechanics UB RAS}, \orgaddress{\city{Yekaterinburg}, \postcode{620108}, \country{Russia}}}

\affil[5]{\orgname{Ural Federal University}, \orgaddress{\city{Yekaterinburg}, \postcode{620002}, \country{Russia}}}

\affil[6]{\orgdiv{School of Science}, \orgname{Jiangnan University}, \orgaddress{\city{Wuxi}, \postcode{214122}, \country{P. R. China}}}


\abstract{If $G$ is a finite group, then the spectrum $\omega(G)$ is the set of all element orders of $G$. The prime spectrum $\pi(G)$ is the set of all primes belonging to $\omega(G)$. A simple graph $\Gamma(G)$ whose vertex set is $\pi(G)$ and in which two distinct vertices $r$ and $s$ are adjacent if and only if $rs \in \omega(G)$ is called the Gruenberg--Kegel graph or the prime graph of $G$.
	
	In this paper, we prove that if $G$ is a group of even order, then the set of vertices which are non-adjacent to $2$ in $\Gamma(G)$ forms a union of cliques. Moreover,  we decide when a strongly regular graph is isomorphic to the Gruenberg--Kegel graph of a finite group.}




\maketitle

\hfill To Professor Otto Kegel on the occasion of his 90th birthday

\section{Introduction}\label{sec1}

Throughout the paper we consider only finite groups and simple graphs, and henceforth the term {\it group} means finite group and the term {\it graph} means simple graph, that is, an undirected graph without loops and multiple edges.

\medskip

If $G$ is a group, then the {\it spectrum} $\omega(G)$ is the set of all element orders of $G$. The {\it prime spectrum} $\pi(G)$ is the set of all primes belonging to $\omega(G)$. A graph $\Gamma(G)$ whose vertex set is $\pi(G)$ and in which two distinct vertices~$r$ and~$s$ are adjacent if and only if $rs \in \omega(G)$ is called the {\it Gruenberg--Kegel graph} or the {\it prime graph} of~$G$. Denote the number of connected components of $\Gamma(G)$
by $s(G)$, and the set of connected
components of $\Gamma(G)$ by $\{\pi_i(G) \mid 1 \leq i \leq s(G) \}$; for a group $G$
of even order, we assume that $2 \in \pi_1(G)$. Denote by $t(G)$ the \emph{independence number} of $\Gamma(G)$, that is, the maximal size of a coclique  (i.\,e. induced subgraph with no edges) in $\Gamma(G)$. If $r\in\pi(G)$, then denote by $t(r,G)$ the maximal size of a coclique in $\Gamma(G)$ containing $r$.

\medskip

Recently the question of characterization of a finite group by its Gruenberg--Kegel graph is under active  investigation. A survey of recent results in this direction can be found, for example, in~\cite{CaMas,Mas_Pans_Star}. For the question of characterization by Gruenberg--Kegel graph, the cases when the graph is connected and when the graph is disconnected are fundamentally different. If the graph is disconnected, then the Gruenberg--Kegel theorem (see Lemma~\ref{Gruenberg--Kegel theorem} below) is a helpful tool. In the case of a connected graph, the situation is more complicated. In this case there exists a strong generalization of the Gruenberg--Kegel theorem proved by A.~Vasil'ev, see Lemma~\ref{vas} below. In this paper, we continue the investigation of the structure of Gruenberg--Kegel graphs of finite groups and prove the following theorem which generalizes Lemma~\ref{vas}.

\medskip

\begin{theorem}\label{Vertices-non-adjacent-to-2} Let $G$ be a finite group of even order such that $t(2,G)\ge 2$. Let $\tau$ be the set of vertices of $\Gamma(G)$ which are not adjacent to $2$. Then the following statements hold{\rm:}
	
	$(1)$ If $G$ is non-solvable, then $G$ has the following normal series $$1 \unlhd K \unlhd G_0 \unlhd G,$$ where $K$ is the largest solvable normal subgroup of $G$, $G_0/K \cong S$ is a finite non-abelian simple group and $G/K$ is almost simple with socle $S$ and either $\tau \subseteq \pi(K) \setminus \pi(G/K)$ or $\tau \subseteq \pi(S)\setminus (\pi(K) \cup \pi(G/G_0))$. In particular, $t(2,G)=2$ or $t(2,G)\le t(2,S)$.
	
	$(2)$ $\tau$ is a union of cliques.
\end{theorem}

\medskip

{\bf Remark~1.} Note that if $p$ is an odd prime, then there exists a finite group such that the set of the vertices which are not adjacent to $p$ in $\Gamma(G)$, is connected and is not a clique. Indeed, let $G=PGL_2(p^m)$, where $m \ge 5$. Then by \cite{Butur}, $\omega(G)$ consists of the divisors of the numbers from the following set $\{p, p^m+1, p^m-1\}$. By the Bang--Zsigmondy Theorem (see Lemma~\ref{zsigm} below), $|\pi(p^m-1)|>1$ and $|\pi(p^m+1)|>1$. Thus, the set $\pi(p^m-1)\cup \pi(p^m+1)$ of primes which are not adjacent to $p$ in $\Gamma(G)$ is connected and is not a clique in $\Gamma(G)$.

\medskip

A number of recent papers are devoted to investigation of combinatorial properties of Gruenberg--Kegel graphs of finite groups. For example, Gruber et.\,al.~\cite{Gruber_et_al} have proved that a graph is isomorphic to the Gruenberg--Kegel graph of a solvable group if and only if its complement is $3$-colorable and triangle-free. In this paper we continue such investigations.

\medskip

A graph $\Gamma$ is called {\it $k$-regular} if each vertex degree of $\Gamma$ is equal to $k$. A {\it strongly regular graph} with parameters $(v, k, \lambda, \mu)$ is a connected $k$-regular graph $\Gamma$ with $v$ vertices such that every two adjacent vertices have $\lambda$ common neighbours and every two non-adjacent vertices have $\mu$ common neighbours for some integers $\lambda \ge 0$ and $\mu \ge 1$. Note that the four parameters $(v, k, \lambda, \mu)$ in a strongly regular graph are not independent. They must obey the following relation:
$$(v-k-1)\mu =k(k-\lambda -1).$$ The complement of a strongly regular graph with parameters $(v, k, \lambda, \mu)$ is either a union of cliques or a strongly regular graph with parameters $$(v, v - k - 1, v - 2 - 2k + \mu, v - 2k + \lambda).$$

The following question was asked by Jack Koolen in a private communication with the third author.

\medskip

\begin{question}[J. Koolen, 2016]\label{Quest_Koolen} What are strongly regular graphs which are isomorphic to Gruenberg--Kegel graphs of finite groups{\rm?}
\end{question}

\medskip

As a corollary of Theorem~\ref{Vertices-non-adjacent-to-2}, we prove the following theorem which gives an answer to Question~\ref{Quest_Koolen}.

\medskip

\begin{theorem}\label{Strongly_regular-graphs}
	Let $\Gamma$ be a strongly regular graph such that $\Gamma$ is isomorphic to the Gruenberg--Kegel graph $\Gamma(G)$ of a finite group $G$. Then one of the following statements holds{\rm:}
	
	$(1)$ $\Gamma$ is the complement to a triangle-free strongly regular graph{\rm;}
	
	$(2)$ $\Gamma$ is a complete multipartite graph with all parts of size $2$.
\end{theorem}

\medskip

{\bf Remark~2.} Note that the complement to a complete multipartite graph with all parts of size $2$ is $2$-colorable and triangle-free. Thus, by~{\rm\cite{Gruber_et_al}}, each graph from Statement~$(2)$ of Theorem~{\rm\ref{Strongly_regular-graphs}} is isomorphic to the Gruenberg--Kegel graph of a solvable group.

\medskip

An important step in proof of Theorem~\ref{Strongly_regular-graphs} is the following theorem.

\medskip

\begin{theorem}\label{ComplMult}
	
	Let $\Gamma$ be a complete multipartite graph with each part of size at least $3$.
	Then $\Gamma$ is not isomorphic to the Gruenberg--Kegel graph of a finite group.
\end{theorem}

\medskip

Note that in \cite{Masl_Pagon}, the third author and D.~Pagon have proved Theorem~\ref{ComplMult} for complete bipartite graphs, i.\,e. they have proved that a complete bipartite graph $K_{n,m}$ is isomorphic to the Gruenberg--Kegel graph of a finite group if and only if $m+n \le 6$ and $(n,m)\not = (3,3)$.

\section{Preliminaries}

Let $n$ be an integer, $\pi$ be a set of primes and $G$ be a group. Denote by $\pi(n)$ the set of all prime divisors of $n$. Note that with respect to this notation, $\pi(G)=\pi(|G|)$. The largest divisor $m$ of $n$ such that $\pi(m)\subseteq \pi$ is called the {\it $\pi$-part} of $n$ and is denoted by $n_\pi$. By $\pi'$ we denote the set of primes which do not belong to $\pi$. If $\pi$ consists of a unique element $p$, then we will write $n_p$ and $n_{p'}$ instead of $n_{\{p\}}$ and $n_{\{p\}'}$, respectively.  $G$ is called a $\pi$-group if $\pi(G) \subseteq \pi$. A subgroup $H$ of $G$ is called a $\pi$-Hall subgroup if $\pi(H)\subseteq \pi$ and $\pi(|G:H|)\subseteq \pi'$.

We will denote by $S(G)$ the {\it solvable radical} of $G$ (the largest solvable normal subgroup of $G$), by $F(G)$ the {\it Fitting subgroup} of $G$ (the largest nilpotent normal subgroup of $G$),
and by $Soc(G)$ the {\it socle} of $G$ (the subgroup of $G$ generated by the set of all non-trivial minimal normal subgroups of $G$), by $O_\pi(G)$ the largest normal $\pi$-subgroup of $G$, and if $\pi=\{p\}$, then we write $O_p(G)$ instead of $O_{\{p\}}(G)$.


\medskip

If $n$ is an integer and $r$ is an odd prime with $(r, n) = 1$, then $e(r, n)$ denotes the {\it multiplicative order} of $n$ modulo $r$. Given an odd integer $n$, we put $e(2, n) = 1$ if $n\equiv1\pmod{4}$, and $e(2,n)=2$ otherwise.

The following lemma is proved in~\cite{Bang}, and also in~\cite{Zs92}.

\medskip

\begin{lem}[{\rm Bang--Zsigmondy}]\label{zsigm} Let $q$ be an integer greater than $1$. For every positive integer $m$ there exists a prime $r$ with $e(r,q)=m$ except in the cases $q=2$ and $m=1$, $q=3$ and $m=1$, and $q=2$ and $m=6$. \end{lem}	

\medskip

Fix an integer $a$ with $|a|>1$. A prime $r$ is said to be a {\it primitive prime divisor} of $a^i-1$ if $e(r,a)=i$. We write $r_i(a)$ (or just $r_i$ if $a$ has been fixed) to denote some primitive prime divisor of $a^i-1$ if such a prime exists, and $R_i(a)$ to denote the set of all such divisors.

\medskip

Let $m$ be a positive integer. Following \cite{VaVd05}, define $$\nu(m)=\begin{cases} m, \mbox{ }m \equiv 0\pmod{4},\\ m/2, \mbox{ }m \equiv 2\pmod{4}, \\
	2m, \mbox{ }m \equiv 1\pmod{2};\end{cases} \eta(m)=\begin{cases} m, \mbox{ }m \equiv 1\pmod{2},\\ m/2, \mbox{ }m \equiv 0\pmod{2}.\end{cases}$$

\medskip

A group $G$ is called a {\em Frobenius group} if there is a subgroup $H$ of $G$ such that $H \cap H^g=1$ whenever $g \in G\setminus H$. Let $$
K=\{1_G\} \cup (G \setminus ( \cup_{g \in G} H^g))$$ be the {\em Frobenius kernel} of $G$. It is well-known (see, for example,~\cite[35.24 and~35.25]{Asch86}) that $K \trianglelefteq G$, $G=K \rtimes H$, $C_G(h)\le H$ for each $h \in H$, and $C_G(k)\le K$ for each $k \in K$.
Moreover, by the Thompson theorem on finite groups with fixed-point-free automorphisms of prime order~\cite[Theorem~1]{Thompson}, $K$ is nilpotent.

A $2$-Frobenius group is a group $G$ which contains a normal Frobenius subgroup $R$ with Frobenius kernel $A$ such that $G/A$ is a Frobenius group with Frobenius kernel $R/A$.

\medskip

\begin{lem}[{\rm Gruenberg--Kegel Theorem, \cite[Theorem~A]{Williams}}]\label{Gruenberg--Kegel theorem} If~$G$ is a group with disconnected Gruenberg--Kegel graph, then one of the following statements holds{\rm:}

		$(1)$ $G$ is a Frobenius group{\rm;}

		$(2)$ $G$ is a $2$-Frobenius group{\rm;}
	
	$(3)$ $G$ is an extension of a nilpotent $\pi_1(G)$-group by a group~$A$, where $S \unlhd A\le\Aut(S)$,~$S$ is a non-abelian simple group with $s(G)\le s(S)$, and $A/S$ is a $\pi_1(G)$-group.
	\
\end{lem}

\medskip

\begin{lem}[{\rm \cite{Va05}}]\label{vas}
	Let $G$ be a non-solvable group with $t(2,G)\geq2$. Then the following statements hold.
	
	$(1)$ There exists a non-abelian simple group $S$ such that $S \unlhd \overline{G} = G/K \le\operatorname{Aut}(S)$,  where $K$ is the solvable radical of $G$.
	
	$(2)$ For every coclique $\rho$ of $\Gamma(G)$ of size at least three, at most one prime in $\rho$ divides the product $|K|\cdot|\overline{G}/S|$. In particular, $t(S)\geq t(G)-1$.
	
	$(3)$ One of the following two conditions holds:	
	
	$\mbox{ }$$\mbox{ }$$\mbox{ }$$(3.1)$ $S\cong Alt_7$ or $L_2(q)$ for some odd $q$, and $t(S)=t(2,S)=3$.
	
	$\mbox{ }$$\mbox{ }$$\mbox{ }$$(3.2)$
	Every prime $p\in\pi(G)$ non-adjacent to $2$ in $\Gamma(G)$ does not divide the product $|K|\cdot|\overline{G}/S|$. In particular, $t(2,S)\geq t(2,G)$.
\end{lem}

\medskip

The following assertion is easy to prove, and can be found, for example, in~\cite[Theorem~1]{Lucido2}.

\medskip

\begin{lem}\label{non-solvable} Let $G$ be a group with ${t(G) \ge 3}$. Then $G$ is non-solvable.
\end{lem}

\medskip

\begin{lem}\label{NormalSeriesAdj} Let $A$ and $B$ be normal subgroups of a group $G$ such that $A\le B$.
	If $r, s \in \pi(B/A)\setminus (\pi(A) \cup \pi(G/B))$, then $r$ and $s$ are adjacent in $\Gamma(G)$ if and only if $r$ and $s$ are adjacent in $\Gamma(B/A)$.
\end{lem}

\medskip

{\it Proof.} The proof of this lemma is elementary.
\medskip

A subgroup $H$ is {\it pronormal} in a group $G$ if the subgroups $H$ and $H^g$ are conjugate in the subgroup $\langle H, H^g \rangle$ for each $g \in G$.

\medskip

\begin{lem}[{\rm \cite[Lemma~4]{GuoMasRev}}]\label{pronormal}
	Let $H \le A$ and $A \unlhd G$. The following statements are equivalent{\rm:}
	
	$(1)$ $H$ is pronormal in $G${\rm;}
	
	$(2)$ $H$ is pronormal in $A$ and $G=AN_G(H)$.
	
\end{lem}

\medskip

\begin{lem}[{\rm \cite[Lemma 1]{Staroletov2}}]\label{semid}
	Let $N$ be an elementary abelian normal subgroup of a group $G$ and $H = G/N$. Define a homomorphism $\phi: H\rightarrow\operatorname{Aut}(N)$ as follows $n^{\phi(gN)}=n^g$. Then $\Gamma(G)=\Gamma(N\rtimes_{\phi}H)$.	
\end{lem}

\medskip

Let $\Gamma$ be a graph, $V(\Gamma)$ be the vertex set of $\Gamma$, and $u \in V(\Gamma)$. Denote by $N(u)$ the set of all vertices which are adjacent to $u$ in $\Gamma$, and by $N_2(u)$ the set of vertices which are at distance $2$ from $u$ in $\Gamma$. It is well-known that each strongly regular graph has diameter $2$. Thus, if $\Gamma$ is strongly regular, then for each $u \in V(\Gamma)$,

$$V(\Gamma)=\{u\}\cup N(u)\cup N_2(u).$$

\medskip

\begin{lem}[{\rm \cite[Lemma~3.1]{SRG_Gardiner_et_al}}]\label{SRG_DSN}
	Let $\Gamma$ be a strongly regular graph. If there exists $u \in V(\Gamma)$ such that $N_2(u)$ is disconnected, then $\Gamma$ is a complete multipartite graph with parts of the same size.
	
\end{lem}

\section{Vertices which are non-adjacent to $2$ in the Gruenberg--Kegel graph of a finite group}

The aim of this section is to prove Theorem~\ref{Vertices-non-adjacent-to-2}. We do this via the following series of assertions.

\smallskip

The following two propositions can be proved by following the arguments in \cite{Va05} and \cite{VaGo09}. Here we provide proofs that are case-free.

\medskip

\begin{prop}\label{2_non_adjacent}
	
	Let $G$ be a group and $K$ be a solvable normal subgroup of $G$ such that $G/K \cong S$ is non-abelian simple. Then in $\Gamma(G)$, $2$ is adjacent to each odd prime from $\pi(K)$ or $2$ is adjacent to each odd prime from $\pi(S)$.
	
\end{prop}

\medskip

{\it Proof.} Let $G$ be a minimal counterexample and let $\tau$ be the set of odd vertices, which are non-adjacent to $2$ in $\Gamma(G)$.

Assume that $r \in \tau \cap \pi(K)$ and $\tau \cap \pi(S) \not = \varnothing$. Let $H$ be a $\{2,r\}$-Hall subgroup of $K$. Take any $g \in G$. Then $H^g \le K$ and the subgroup $\langle H, H^g \rangle \le K$ is solvable, therefore by the Hall theorem, $H$ and $H^g$ are conjugate in $\langle H, H^g \rangle$. Thus, $H$ is pronormal in $G$ and by Lemma~\ref{pronormal}, $G=KN_G(H)$. We have $G/K\cong N_G(H)/N_K(H)$, therefore $\pi(S) \cup \{2,r\} \subseteq \pi(N_G(H))$ and $\Gamma(N_G(H))$ is a subgraph of $\Gamma(G)$. Thus, $N_G(H)$ is a counterexample to the statement of the proposition, therefore $N_G(H)=G$ by minimality of $G$.

We now show that $|H|$ is even. If $|H|$ is odd, then each Sylow $2$-subgroup of $G$ is isomorphic to a Sylow $2$-subgroup of $G/K$, therefore by the Glauberman $Z^*$-theorem~\cite{Glauberman}, $G$ has a subgroup isomorphic to the Klein $4$-group, therefore by \cite[Theorem~10.3.1]{Gorenstein}, a Sylow $2$-subgroup of $G$ can not act fixed-point-freely on any group of odd order; a contradiction. Now by Lemma~\ref{Gruenberg--Kegel theorem}, we have that $H$ is either a Frobenius group or a $2$-Frobenius group and $\pi(H)=\{2, r\}$.

Suppose that $H$ is a Frobenius group with Frobenius kernel $A$, which is a $2$-group and Frobenius complement of odd order. Then $G/A$ is also a counterexample to the proposition and $|G/A|<|G|$, contradicting to the minimality of $G$.

Suppose that $H$ is a $2$-Frobenius group, where $R$ is a normal subgroup of $K$, with the property that $R$ is a Frobenius group with Frobenius kernel $A$. Then $G/A$ is also a counterexample to the proposition and $|G/A|<|G|$, contradicting to the minimality of $G$.

Thus, by minimality of $G$, $H$ is a Frobenius group with Frobenius kernel $F=F(H)$ such that $\pi(F)=\{r\}$ and Frobenius complement $D$ such that $\pi(D)=\{2\}$.  Also by minimality of $G$ we can assume that $F$ is an elementary abelian $r$-group, therefore by Lemma~\ref{semid}, we can assume that $G=F \rtimes C$, where $C \cong G/F$ and in $\Gamma(C)$, $2$ is adjacent to each odd prime in $\pi(S(C))$. Moreover, $H/F\cong O_2(C)$ and $|S(C)/O_2(C)|$ is odd.

Let $Q$ be a Sylow $2$-subgroup of $C$. Since $Q$ acts on $F$ fixed-point-freely by \cite[Theorem~10.3.1]{Gorenstein}, we have $Q$ is either cyclic or generalized quaternion. In any case, $Q$ has a unique involution $i$ which is contained in each normal subgroup of $Q$, in particular, $i$ is contained in $Z(Q)$ and in $S(C)\cap Q=O_2(C)$. Thus, $i$ is contained in the center of each conjugate of $Q$. Consider a subgroup $W=\langle Q^c\mid c \in C\rangle$
of $C$ which is normal in $C$ and such that $|C:W|$ is odd. Then by the Feit-Thompson theorem, $C/W$ is solvable. Thus, $W$ contains a unique non-abelian composition factor of $C$ and therefore $\pi(S)\subseteq \pi(W)$. But $i$ is contained in $Z(W)$, therefore, $2$ is adjacent to each odd prime from $\pi(W)$ in $\Gamma(W)$ and so, $2$ is adjacent to each odd prime from $\pi(S)$ in $\Gamma(G)$, contradicting the assumption that $G$ is a counterexample to the statement of the proposition. \hfill $\Box$

\medskip

\begin{prop}\label{anti2}
	Let $G$ be a non-solvable group and $\tau$ be the set of vertices which are not adjacent to $2$ in $\Gamma(G)$. If $|\tau|\ge 1$, then $G$ has the following normal series $$1 \unlhd K \unlhd G_0 \unlhd G,$$ where $K=S(G)$ is the solvable radical of $G$, $G_0/K \cong S$ is a non-abelian simple group and $G/K$ is almost simple with socle $S$ such that  either $\tau \subseteq \pi(K) \setminus \pi(G/K)$ or $\tau \subseteq \pi(S)\setminus (\pi(K) \cup \pi(G/G_0))$.
\end{prop}

\medskip

{\it Proof.} By Lemma~\ref{vas}, if $G$ is non-solvable and $|\tau|\ge 1$, then $G$ has the following normal series $$1 \unlhd K \unlhd G_0 \unlhd G,$$ where $K$ is the solvable radical of $G$, $G_0/K \cong S$ is a non-abelian simple group and $G/K$ is almost simple with socle $S$.

Let $t \in \tau$. It is clear that $\Gamma(G/K)$ is a subgraph of $\Gamma(G)$, therefore by \cite[Lemma~1.2]{Va05}, we have ${\tau \cap \pi(G/G_0) = \varnothing}$, and therefore we have $$\tau \subseteq (\pi(S) \cup \pi(K)) \setminus \pi(G/G_0).$$

By Proposition~\ref{2_non_adjacent}, we have $\tau \cap \pi(K) = \varnothing$ or $\tau \cap \pi(S)=\varnothing$. Thus, either $\tau \subseteq \pi(K)\setminus \pi(G/K)$ or $\tau \subseteq \pi(S)\setminus (\pi(K) \cup \pi(G/G_0))$. \hfill $\Box$

\medskip

\begin{lem}\label{2_non-adj_simple} Let $G$ be a non-abelian simple group and $\tau$ be the set of vertices which are not adjacent to $2$ in $\Gamma(G)$. If $\tau \not = \varnothing$, then $\tau$ is a union of cliques.
\end{lem}

\medskip

{\it Proof.} Proof of the lemma for sporadic simple groups follows directly from~\cite{Atlas}.

\smallskip

Assume that $G=Alt_n$ for $n \ge 5$ and $p$ is non-adjacent to $2$ in $\Gamma(G)$. Then $p+4 > n$ and therefore $n \ge p > n-4$. Note that between the numbers $n$, $n-1$, $n-2$ and $n-3$ there are at most $2$ odd primes. Thus, if $\tau \not = \varnothing$, then $\tau$ consists either from an only prime or from exactly two primes and is a union of cliques in any case.

\smallskip

Let $G$ be a group of Lie type. Consider the possibilities for $G$ case by case.

If $G=A_n(q)$, then the statement of lemma follows directly from the adjacency criterion for $\Gamma(G)$, which can be found in \cite[Propositions~2.1,~3.1, and~4.1]{VaVd05}.

If $G={^2}A_n(q)$, then the statement of lemma follows directly from the adjacency criterion for $\Gamma(G)$, which can be found in \cite[Propositions~2.2,~3.1, and~4.2]{VaVd05}.

If $G=B_n(q)$ or $G=C_n(q)$, then the statement of lemma follows directly from the adjacency criterion for $\Gamma(G)$, which can be found in \cite[Propositions~3.1 and~4.3]{VaVd05} and \cite[Proposition~2.4]{VaVd11}.

If $G=D_n(q)$ or $G={^2}D_n(q)$, then the statement of lemma follows directly from the adjacency criterion for $\Gamma(G)$, which can be found in \cite[Propositions~3.1 and~4.4]{VaVd05} and \cite[Proposition~2.5]{VaVd11}.

If $G$ is isomorphic to one of the groups $E_8(q)$, $E_7(q)$, $E_6(q)$, ${^2}E_6(q)$,
$F_4(q)$, ${^3}D_4(q)$, $G_2(q)$, then the statement of lemma follows directly from the adjacency criterion for $\Gamma(G)$, which can be found in \cite[Propositions~3.2 and~4.5]{VaVd05} and \cite[Proposition~2.7]{VaVd11}.

If $G$ is isomorphic to one of the groups ${^2}F_4(q)'$, ${^2}G_2(q)$, $^{2}B_2(q)$, then the statement of lemma follows directly from the adjacency criterion for $\Gamma(G)$, which can be found in \cite[Propositions~3.3 and~4.5]{VaVd05} and \cite[Proposition~2.9]{VaVd11}. \hfill $\Box$

\medskip

\begin{prop}\label{cocliques} Let $G$ be a group of even order and $\tau$ be the set of vertices which are not adjacent to $2$ in $\Gamma(G)$. If $\tau \not = \varnothing$, then $\tau$ is a union of cliques.
\end{prop}

\medskip

{\it Proof.} Assume that $H$ is solvable.  Let $\tau_1=\tau \cup \{2\}$ and $H$ be a $\tau_1$-Hall subgroup of $G$. Then by the Hall theorem, primes $p$ and $q$ from $\tau_1$ are adjacent in $\Gamma(G)$ if and only if they are adjacent in $\Gamma(H)$. Note that $H$ is solvable and $2$ is an isolated vertex in $\Gamma(H)$. Now if $\tau$ does not form a clique in $\Gamma(H)$, then Lemma~\ref{non-solvable} gives a contradiction with solvability of $H$.

Thus, we can assume that $G$ is non-solvable and $\tau \not = \varnothing$. By Lemma~\ref{vas}, $$S \unlhd \overline{G} = G/K \le\operatorname{Aut}(S),$$  where $K$ is solvable and $S$ in a non-abelian simple group. By Proposition~\ref{anti2},  either $\tau \subseteq \pi(K) \setminus \pi(G/K)$ or $\tau \subseteq \pi(S)\setminus (\pi(K) \cup \pi(\overline{G}/S))$.

Let $\tau \subseteq \pi(K)$. Consider a Sylow $2$-subgroup $Q$ of $G$ and let $H=KQ$. It is clear that $H$ is solvable and $\Gamma(H)$ is a subgraph in $\Gamma(G)$. From above, $\tau$ forms a clique in $\Gamma(H)$, therefore $\tau$ is a clique in $\Gamma(G)$.

Let $\tau \subseteq \pi(S)\setminus (\pi(K) \cup \pi(\overline{G}/S))$. By Lemma~\ref{NormalSeriesAdj}, the primes $p, q\in \tau$ are non-adjacent in $\Gamma(G)$ if and only if they are non-adjacent in $\Gamma(S)$. Moreover, $\tau$ is a subset of the set $\sigma$ of odd vertices which are non-adjacent to $2$ in $\Gamma(S)$. By Lemma~\ref{2_non-adj_simple}, $\sigma$ is a union of cliques. Thus, $\tau$ is an induced subgraph of a union of cliques, therefore, $\tau$ is also a union of cliques. $\phantom{x}$ \hfill $\Box$

\medskip

Theorem~\ref{Vertices-non-adjacent-to-2} follows directly from Propositions~\ref{anti2} and~\ref{cocliques}. \hfill $\Box$

\medskip

\section{Strongly regular graphs which are isomorphic to Gruenberg--Kegel graphs of finite groups}

The aim of this section is to prove Theorems~\ref{Strongly_regular-graphs} and~\ref{ComplMult}.

\medskip

{\it Proof of Theorem~{\rm\ref{Strongly_regular-graphs}}.} Let $\Gamma$ be a strongly regular graph and suppose that there exists a group $G$ such that $\Gamma$ is isomorphic to $\Gamma(G)$. If $|G|$ is odd, then by the Feit--Thompson theorem, $G$ is solvable. Thus, by Lemma~\ref{non-solvable}
the complement of $\Gamma$ is either a triangle-free strongly regular graph or a union of cliques of size $2$. So, we can assume that $|G|$ is even. Then by Theorem~\ref{Vertices-non-adjacent-to-2}, $N_2(2)$ is a union of cliques.

\medskip

Assume that $N_2(2)$ is a clique. Then from strong regularity of $\Gamma$, it follows that $N_2(v)$ is a clique for each vertex $v$ of $\Gamma$. Thus, in the complement of $\Gamma$ any two adjacent vertices do not have a common neighbour. Therefore, again $\Gamma$ is a complete multipartite graph with parts of size $2$ or the complement of $\Gamma$ is a triangle-free strongly regular graph.

\medskip

Assume that $N_2(2)$ is a union of more than one cliques. Then by Lemma~\ref{SRG_DSN}, $\Gamma$ is a complete multipartite graph with parts of the same size $t$.

\medskip

Now to complete proof of Theorem~\ref{Strongly_regular-graphs} it is sufficient to prove Theorem~\ref{ComplMult}.

\medskip

{\it Proof of Theorem~{\rm\ref{ComplMult}}.} Let $\Gamma$ be a complete multipartite graph with each part of size at least $3$. Assume that $\Gamma$ is isomorphic to $\Gamma(G)$, where $G$ is a group. Then $t(G) \ge 3$ and therefore by Lemma~\ref{non-solvable}, $G$ is non-solvable. By the Feit--Thompson theorem, $|G|$ is even.

Let $\sigma$ be a part of $\Gamma(G)$ which contains $2$, and then $\tau =\sigma \setminus \{2\}$ is exactly the set of the vertices which are not adjacent to $2$ in $\Gamma$. By Theorem~\ref{Vertices-non-adjacent-to-2}, $G$ has the following normal series $$1 \unlhd K \unlhd G_0 \unlhd G,$$ where $K$ is solvable, $G_0/K \cong S$ is a non-abelian simple group and $G/K$ is almost simple with socle $S$ and $\tau \subseteq \pi(S)\setminus (\pi(K) \cup \pi(G/G_0))$. In particular, by Lemma~\ref{NormalSeriesAdj}, $\sigma$ forms in $\Gamma(S)$ a coclique containing the vertex $2$ and $t(2,S) \ge t(2,G)=|\sigma|\ge 3$.

Let $\mu\not =\sigma$ be another part of $\Gamma$. Since $\mu$ is a coclique, by Statement~$(2)$ of Lemma~\ref{vas}, there are at least two primes $x$ and $y$ with $$\{x, y\}\subseteq \mu \cap (\pi(S) \setminus (\pi(K) \cup \pi(G/G_0))).$$ Therefore $x$ and $y$ are non-adjacent in $\Gamma(G)$ and by Lemma~\ref{NormalSeriesAdj}, each vertex from $\{x, y\}$ is adjacent in $\Gamma(S)$ to each vertex from $\tau$.

Assume that $\Gamma(S)$ is disconnected. Let $i>1$ and  $u \in \tau \cap \pi_i(S)$. If $x$ is adjacent to $u$ and $y$ is adjacent to $u$ in $\Gamma(S)$, we have $\{x,y\} \subseteq \pi_i(S)$ and therefore $x$ and $y$ are adjacent, a contradiction. Thus, $\tau \cap \pi_i(S) = \varnothing$ for each $i>1$. Similarly, $\{x,y\} \cap \pi_i(S) = \varnothing$ for each $i>1$.

\smallskip

By \cite{Atlas}, $S$ is not a sporadic simple group, and it is clear that $S$ is not an alternating simple group. Thus, $S$ is a simple group of Lie type. From \cite[Tables~1--3]{AlexKond}, \cite[Propositions~2.1, 2.2, 3.1, 3.2, 3.3, 4.1, 4.2, 4.3, 4.4, and~4.5]{VaVd05}, \cite[Proposition~2.4, 2.5, 2.7, and~2.9]{VaVd11} and \cite[Tables~2--7]{VaVd05}, taking into account corrections from \cite[Appendix]{VaVd11}, we conclude that one of the following statements holds{\rm:}

\begin{itemize}
	
	\item[(1)] $S \cong A_1(q)$, or $q$ is even and $S \cong A_2(q)$ or ${^2}A_2(q)$ with some extra conditions on $q${\rm;}
	
	\item[(2)]  $S\cong A_{n-1}(q)$, $n>3$, $q$ is even, $(n,q)\not =(6,2), (7,2)$, and $\tau = \{r_{n-1},r_n\}${\rm;}
	
	\item[(3)]  $S\cong {^2}A_{n-1}(q)$, $n>3$, $q$ is even, and one of the following statements holds{\rm:}
	
	\begin{itemize}
		
		\item[(3i)] $n \equiv 0 \pmod{4}$, $(n,q)\not = (4,2)$, and $\tau=\{r_{2n-2}, r_n\}${\rm;}
		
		\item[(3ii)] $n \equiv 1 \pmod{4}$, and $\tau=\{r_{n-1}, r_{2n}\}${\rm;}
		
		\item[(3iii)] $n \equiv 2 \pmod{4}$, and $\tau=\{r_{2n-2}, r_{n/2}\}${\rm;}
		
		\item[(3iv)] $n \equiv 3 \pmod{4}$, and $\tau=\{r_{(n-1)/2}, r_{2n}\}${\rm;}
		
	\end{itemize}
	
	\item[(4)]  $S\cong B_n(q)$ or $C_n(q)$, $n>1$ is odd, $q$ is even, $(n,q)\not = (3,2)$, and $\tau=\{r_n, r_{2n}\}${\rm;}

	\item[(5)]  $S\cong D_n(q)$, $n\ge 4$, $q$ is even, $(n,q)\not = (4,2)$, and one of the following statements holds{\rm:}
	
	\begin{itemize}
		
		\item[(5i)] $n \equiv 0 \pmod{2}$ and $\tau=\{r_{n-1}, r_{2n-2}\}${\rm;}
		
		\item[(5ii)] $n \equiv 1 \pmod{2}$, and $\tau=\{r_n, r_{2n-2}\}${\rm;}

	\end{itemize}

	\item[(6)]  $S\cong {^2}D_n(q)$, $n\ge 4$, $q$ is even, $(n,q)\not = (4,2)$, and one of the following statements holds{\rm:}
	
	\begin{itemize}
		
		\item[(6i)] $n \equiv 0 \pmod{2}$ and $\tau \subseteq\{r_{n-1}, r_{2n-2}, r_{2n}\}${\rm;}
		
		\item[(6ii)] $n \equiv 1 \pmod{2}$, and $\tau=\{ r_{2n-2}, r_{2n}\}${\rm;}
		
	\end{itemize}

	\item[(7)]  $S\cong E_7(q)$, $q$ is even, and $\tau \subseteq \{r_7, r_9, r_{14}, r_{18}\}${\rm;}

	\item[(8)]  $S\cong E_7(q)$, $q$ is odd, and one of the following statements holds{\rm:}
	
	\begin{itemize}
		
		\item[(8i)] $q \equiv 1 \pmod{4}$ and $\tau=\{r_{14}, r_{18}\}${\rm;}
		
		\item[(8ii)] $q \equiv 3 \pmod{4}$, and $\tau=\{r_7, r_9\}${\rm;}

	\end{itemize}

	\item[(9)]  $S\cong A_{n-1}(q)$, $n>3$, $q$ is odd, $n_2=(q-1)_2>2$, and $\tau = \{r_{n-1},r_n\}${\rm;}
	
	\item[(10)]  $S\cong {^2}A_{n-1}(q)$, $n>3$, $q$ is odd, $n_2=(q+1)_2>2$, and $\tau = \{r_{2n-2},r_n\}${\rm;}
	
	\item[(11)]  $S\cong D_n(q)$, $n>4$ is odd, $q \equiv 5 \pmod{8}$, and $\tau = \{r_n,r_{2n-2}\}${\rm;}
	
	\item[(12)]  $S\cong {^2}D_n(q)$, $n>4$ is odd, $q \equiv 3 \pmod{8}$, and $\tau = \{r_{2n-2},r_{2n}\}$.

\end{itemize}

\medskip

Assume that Statement~$(1)$ holds. If $S\cong A_1(q)$, then by \cite[Propositions~2.1,~3.1, and~4.1]{VaVd05}, $\Gamma(S)$ is a union of $3$ cliques. In particular, $\Gamma(S)$ does not contain induced $4$-cycles, a contradiction. If $q$ is even and $S \cong A_2(q)$ or $^2A_2(q)$, then the picture of a compact form of $\Gamma(S)$ can be found on Pic.~7 in \cite{VaVd11}. It is clear that $\Gamma(S)$ does not contain an induced $4$-cycle with two non-adjacent vertices both non-adjacent to $2$, a contradiction.

\medskip

Further let $p$ be the characteristic of the field over which $S$ is defined, $x \in \pi(S)$ with $x$ adjacent to each element from $\tau$, and if $x\not = p$, then put $k=e(x,q)$.

\medskip

Assume that Statement~$(2)$ or Statement~$(9)$ holds. If $p$ is odd, then since $x$ is adjacent both to $r_n$ and $r_{n-1}$ by \cite[Proposition~3.1]{VaVd05}, we have $x\not = p$. Assume that $k>1$. Since $x$ is adjacent to $r_n$, it follows from \cite[Proposition~2.1]{VaVd05} that $k$ divides $n$. On the other hand, $x$ is adjacent to $r_{n-1}$, therefore $k$ divides $n-1$ which is coprime to $n$, a contradiction. Thus, $k=1$. Since $x$ was chosen arbitrarily, each two elements adjacent to both to $r_n$ and $r_{n-1}$ are forced to be adjacent in $\Gamma(S)$, a contradiction.

\medskip

Assume that Statement~$(3)$ or Statement~$(10)$ holds. If $p$ is odd, then since $x$ is adjacent both to $r_n$ (with $\nu(n)=n$) and to $r_{2n-2}$ (with $\nu(2n-2)=n-1$) by \cite[Proposition~3.1]{VaVd05}, $x\not = p$. Using \cite[Proposition~3.1]{VaVd05} if $p=2$ and \cite[Proposition~4.2]{VaVd05} if $p$ is odd, we conclude that $\tau$ consists of two primes $\{a, b\}$ with $\nu(e(a,q))=n$, $\nu(e(b,q))=n-1$. Let $\nu(k)>1$, i.\,e. $k \not =2$. Since $x$ is adjacent to $a$, it follows from \cite[Proposition~2.2]{VaVd05} that $\nu(k)$ divides $n$. On the other hand, $x$ is adjacent to $b$, therefore $\nu(k)$ divides $n-1$ which is coprime to $n$, a contradiction. Thus, $k=2$. Since $x$ was chosen arbitrarily, each two elements adjacent to both to $a$ and $b$ are forced to be adjacent in $\Gamma(S)$, a contradiction.

\medskip

Assume that Statement~$(4)$ holds. Note that $k \le 2n$ and $\eta(k)\le n$. Thus, by \cite[Proposition~2.4]{VaVd11}, if $x$ is adjacent both to $r_n$ and to $r_{2n}$, then $n/k$ and $2n/k$ are both odd integers, a contradiction.

\medskip

Assume that Statement~$(5)$ holds. Up to the end of this paragraph, we refer to \cite[Proposition~2.5]{VaVd11}. If $n$ is odd, then $x$ is adjacent to $r_n$ (with $\eta(n)=n$) and to $r_{2n-2}$ (with $\eta(2n-2)=n-1$). Since  $x$ is adjacent to $r_n$, we have that $n/k$ is an odd integer or $k/n$ is an odd integer, therefore $k$ is odd. Since  $x$ is adjacent to $r_{2n-2}$, we have
$$2(n-1)+2k> 2n - (1-(-1)^{2n-2+k})=2n-2,$$ therefore again $(2n-2)/k$ is an odd integer or $k/(2n-2)$ is an odd integer, a contradiction. If $n$ is even, then $x$ is adjacent to $r_{n-1}$ and to $r_{2n-2}$. We have $$n-1=\eta(n-1)=\eta(2n-2).$$
Moreover, if $k$ is odd, then $$2(n-1)+2k> 2n - (1-(-1)^{2n-2+k}),$$ therefore $(2n-2)/k$ is an odd integer, a contradiction.
If $k$ is even, then $$2(n-1)+k> 2n - (1-(-1)^{n-1+k}),$$ therefore $(n-1)/k$ is an odd integer or $k/(n-1)$ is an odd integer, a contradiction.

\medskip

Assume that Statement~$(11)$ holds. Since $n$ is odd and $x$ is adjacent to $r_n$ (with $\eta(n)=n$) and to $r_{2n-2}$ (with $\eta(2n-2)=n-1$), by \cite[Proposition~3.1]{VaVd05}, $x\not = p$. Then by \cite[Proposition~2.5]{VaVd11}, $n/k$ is an odd integer or $k/n$ is an odd integer, therefore $k$ is odd. Then $$2(n-1)+2k> 2n - (1-(-1)^{2n-2+k})=2n-2,$$ therefore $(2n-2)/k$ is an odd integer or $k/(2n-2)$ is an odd integer, a contradiction.

\medskip

Assume that Statement~$(6)$ holds. Up to the end of this paragraph, we again refer to \cite[Proposition~2.5]{VaVd11}. If $n>4$ is odd, then $x$ adjacent to $r_{2n}$ (with $\eta(2n)=n$) and to $r_{2n-2}$ (with $\eta(2n-2)=n-1$). Since  $x$ adjacent to $r_{2n}$, we have $2n/k$ is an odd integer or $k/2n$ is an odd integer, therefore $k$ is even. Since  $x$ adjacent to $r_{2n-2}$, we have
$$2(n-1)+k> 2n - (1+(-1)^{2n-2+k})=2n-2,$$ therefore again $(2n-2)/k$ is an odd integer or $k/(2n-2)$ is an odd integer, a contradiction. If $n \ge 4$ is even, then $x$ is adjacent to at least two primes from the set $\{r_{n-1}, r_{2n-2}, r_{2n}\}$.  Let $\eta(k)>1$, i.\,e. $k \not \in \{1,2\}.$ Similarly as in Statement~$(5)$ we prove that $r$ can not be adjacent to both $r_{n-1}$ and $r_{2n-2}$. Thus, $x$ is adjacent to $r_{2n}$. This implies that $2n/k$ or $k/2n$ is an odd integer and therefore $k$ is even. If $x$ is adjacent to $r_{2n-2}$, then $$2(n-1)+ k> 2n - (1+(-1)^{2n-2+k})=2n-2,$$ and so, $(2n-2)/k$ is an odd integer or $k/(2n-2)$ is an odd integer, a contradiction. If $x$ is adjacent to $r_{n-1}$, then $$2(n-1)+ k> 2n - (1+(-1)^{n-1+k})=2n,$$ and so, $(n-1)/k$ is an odd integer or $k/(n-1)$ is an odd integer, a contradiction. Thus, $\eta(k)=1$. Since $x$ was chosen arbitrarily, every two distinct vertices adjacent to at least two primes from the set $\{r_{n-1}, r_{2n-2}, r_{2n}\}$ are forced to be adjacent in $\Gamma(S)$, a contradiction.

\medskip

Assume that Statement~$(12)$ holds. Since $n$ is odd and $x$ is adjacent to $r_{2n}$ (with $\eta(n)=n$) and to $r_{2n-2}$ (with $\eta(2n-2)=n-1$), by \cite[Proposition~3.1]{VaVd05}, $x\not = p$. Then by \cite[Proposition~2.5]{VaVd11}, $2n/k$ is an odd integer or $k/2n$ is an odd integer, therefore $k$ is even. Then $$2(n-1)+k> 2n - (1+(-1)^{2n-2+k})=2n-2,$$ therefore $(2n-2)/k$ is an odd integer or $k/(2n-2)$ is an odd integer, a contradiction.

\medskip

Assume that Statement~$(7)$ or Statement~$(8)$ holds. The picture of a compact form of $\Gamma(S)$ can be found on Pic.~4 in \cite{VaVd11}.
It is clear that if $q$ is even, then there are no two non-adjacent vertices in $\Gamma(S)$ that are adjacent to at least two vertices from the set $\{r_7, r_9, r_{14}, r_{18}\}$,  a contradiction. If $q$ is odd, then any two vertices that are adjacent both to $r_7$ and to $r_9$ or both to $r_{14}$ and to $r_{18}$, are adjacent; a contradiction.

\medskip

The proof of Theorem~\ref{ComplMult} is complete.  \hfill $\Box$

\medskip

The proof of Theorem~\ref{Strongly_regular-graphs} is also complete. \hfill $\Box$

\medskip


\section{Acknowledgements}

The authors of this paper contributed equally to this work. The authors are ordered with respect to alphabet ordering in English.

\smallskip

The work was supported by Russian Science Foundation (project no. 24-41-10004, https://rscf.ru/en/project/24-41-10004/ , for example, Theorem~\ref{Vertices-non-adjacent-to-2}), by the Ministry of Science and Higher Education of the Russian Federation, project no. 075-02-2024-1428 for the development of the regional scientific and educational mathematical center ''Ural Mathematical Center'' (for example, Theorem~\ref{Strongly_regular-graphs}), the National Natural Science Foundation of China (project no. 12101166) and the Hainan Provincial Natural Science Foundation of China (project no. 123MS005).

\smallskip

The authors are very thankful to Dr. Alexey Staroletov, Prof. Anotoly Kondrat'ev, and anonymous reviewers for their helpful comments which improved this text.

\end{document}